\documentclass[12pt]{article}
\usepackage{amsmath, amsfonts, amsthm,amssymb, color, hyperref, enumerate, extarrows, cite, subfiles, csquotes, enumitem, commath, mathtools,eucal,mathrsfs, bm }

 \usepackage[titletoc]{appendix}

\newtheorem{Lemma}{Lemma}[section]
\newtheorem{Theorem}[Lemma]{Theorem}
\newtheorem{Proposition}[Lemma]{Proposition}
\newtheorem{Corollary}[Lemma]{Corollary}

\newtheorem{remark}[Lemma]{Remark}
\newtheorem{definition}[Lemma]{Definition}
\newtheorem{example}[Lemma]{Example}

\newtheorem{Fact}[Lemma]{Fact}

\hypersetup{
	colorlinks   = true,
	citecolor    = magenta,
    urlcolor    =blue}

\usepackage[titletoc]{appendix}

	\newcommand{\C}{\mathbb{C}}

 \def\bt{\begin{Theorem}}
\def\et{\end{Theorem}}
\def\bl{\begin{Lemma}}
\def\el{\end{Lemma}}
\def\bp{\begin{Proposition}}
\def\ep{\end{Proposition}}
\def\bcor{\begin{Corollary}}
\def\ecor{\end{Corollary}}
\def\bpf{\begin{proof}}
\def\epf{\end{proof}}

\def\brem{\begin{remark}}
\def\erem{\end{remark}}

\def\bedef{\begin{definition}\rm }
\def\endef{\end{definition}}

\def\beg{\begin{example}}
\def\eeg{\end{example}}

\def\bef{\begin{Fact}}
\def\eef{\end{Fact}}

\def\bc{\begin{center}}
\def\ec{\end{center}}

\def\beq{\begin{equation}}
\def\eeq{\end{equation}}
\def\beqarray{\begin{eqnarray*}}
\def\eeqarray{\end{eqnarray*}}
\def\<{\leftangle}
\def\>{\rightangle}
\def\({\left(}
\def\){\right)}

\def\<{\langle}
\def\>{\rangle}

\def\a{\alpha}
\def\b{\beta}

\def\th{\theta}

\def\e{\varepsilon}

\def\O{\Omega}

\def\z{\zeta}

\def\w.r.t.{with respect to}
\def\R{{\mathbb{R}}}

\def\Z{{\mathbb{Z}}}

\def\B{{\mathbb{B}}}

\def\C{{\mathbb{C}}}

\def\F{{\mathcal{F}}}
\def\A{{\mathcal{A}}}
\def\L{{\mathcal{L}}}

\def\bq{\begin{quote}}
\def\eq{\end{quote}}

\def\bit{\begin{itemize}}
\def\eit{\end{itemize}}

\def\ben{\begin{enumerate}}
\def\een{\end{enumerate}}

\begin{document}

\title{\vspace{-1.2cm} \bf Duality and Interpolation of Bergman Spaces\rm}

\author{Shreedhar Bhat}
\date{}

\maketitle

\begin{abstract}
 This paper explores the dual space corresponding to $p$-Bergman space and examines the essential condition for the dual space to be a $q$-Bergman space. The investigation involves a detailed examination of the interpolation space of a Banach couple. Additionally, we draw comparisons between the `duality', `integrability' and `regularity' properties of a domain.

\end{abstract}

\renewcommand{\thefootnote}{\fnsymbol{footnote}}
\footnotetext{\hspace*{-7mm} 
\begin{tabular}{@{}r@{}p{16.5cm}@{}}
& Keywords. Dual of Bergman Space, Bergman kernel, Interpolation of Banach Spaces.\\
& Mathematics Subject Classification. Primary 32A36; Secondary 32C37; 46B70\end{tabular}}

\section{Introduction} 
For a bounded domain $\Omega\subset \C^n$, $1\leq p<\infty $, define the $p$-Bergman space as  $$\A^p(\O)=\{f\in \mathcal{L}^p(\O):f \text{ is holomorphic in } \O\}$$
For $p=2$, $\A^2(\O)$ is a subspace of the Hilbert space $\L^2(\O)$ and the evaluation functional $e_z(f)=f(z)$ is a continous linear functional. Thus by Riesz representation theorem,  there exists a function $K_z(\cdot) \in \A^2(\O)$ such that 
$$f(z)=\<f, K_z\>$$
The function $K(w,z):=K_z(w)$ is called the \textbf{Bergman kernel} function. Further, the orthogonal projection, called the \textbf{Bergman projection}, from $\L^2(\O)$ to $\A^2(\O)$ is given by $$B: \L^2(\O)\rightarrow \A^2(\O)$$ $$(Bf)(z)=\int_\O f(\z)\cdot K(\z,z)dV_\z.$$

Recently, there has been a growing interest in the study of $p$-Bergman theory. One of the most interesting use of $p$-Bergman space is its role in the classification of bounded hyperconvex domains cf. \cite[Theorem 1.2]{deng2020linear}. Numerous efforts have been undertaken to deepen the comprehension of $p$-Bergman theory, see \cite{chakrabarti2019duality}, \cite{chen2022p} etc.\\

Here we try to investigate the $p$-Bergman space via its dual space. Recall that the dual space of $\L^p(\O)$,  $[\L^p(\O)]'=\L^q(\O)$, when $1/p+1/q=1$. The interesting and the most obvious question would be under what conditions does this duality translate to the $p$-Bergman space i.e. when is $[\A^p(\O)]'=\A^q(\O)$? This question has been studied earlier in \cite{Hedenmalm}, \cite{chakrabarti2019duality}. Before we proceed further, we attempt to make this notion of duality more precise. \\
\newpage 
When $1/p+1/q=1$, define the conjugate-linear map $$\Phi_p: \A^q(\O) \rightarrow A^p(\O)'$$ by $$\Phi_p(g)(f)=\int_\O f\overline{g}= \<f,g\> \text{ for } g\in \A^q(\O), f\in \A^p(\O).$$
Similarly define $$\Phi_q:\A^p(\O)\rightarrow \A^q(\O)'.$$

Thus the duality of $p$-Bergman space relates to the surjectivity of the above function $\Phi_p$. The sufficient conditions for surjectivity are discussed in \cite[Theorem 2.15]{chakrabarti2019duality}. We will look into necessary conditions and therefore assume that both $\Phi_p$ and $\Phi_q$ are surjective throughout the paper.

\bl \label{Isomorphic}
Let $\O$ be a bounded domain in $\C^n$ and $p,q$ be conjugate exponents of each other. Assume that $\Phi_p, \Phi_q$ are surjective. Then $\A^p(\O)$ and $\A^q(\O)'$ are isomorphic to each other i.e., $\Phi_p, \Phi_q$ are injective maps. Similarly, $\A^q(\O)$ and $\A^p(\O)'$ are isomorphic to each other.
\el 
\bpf Without loss of generality, assume that $p>2$ and $q<2$. Since $\Phi_p$ is surjective, it suffices to show that $\Phi_p$ is an injective operator.\\

For $z\in \O$, observe that the evaluation functional $e_z(g)=g(z)$ $\in \A^q(\O)'$. By the assumption of surjectivity, there exists a function $R_q(\cdot,z) \in \A^p(\O)$ such that $$g(z)=\<g,R_q(\cdot,z)\> \text{ for } g\in \A^q(\O)$$
Similarly, there exists a function $R_p(\cdot,z) \in \A^q(\O)$ such that 
$$f(z)=\<f, R_p(\cdot,z)\> \text{ for } f\in \A^p(\O)$$
Thus for $z,w \in \O$,
$$R_p(z,w)=\<R_p(\cdot,w), R_q(\cdot,z))\>=\overline{\<R_q(\cdot,z)), R_p(\cdot,w)\>}=\overline{R_q(w,z)}$$
This shows that there is a unique such reproducing function $R_q(\cdot,*)$ and $R_p(\cdot,*)$ which shows that the linear operators $\Phi_p$ and $\Phi_q$ are injective, thereby proving Lemma \ref{Isomorphic}.
\epf 

\brem \label{Integrability of BK}
Since $p>2$ and $\O$ is a bounded domain in $\C^n$, $\A^p(\O)\subset \A^2(\O)\subset \A^q(\O)$. By the reproducing property of the Bergman kernel, we have $$f(z)=\<f, K(\cdot,z)\> $$ for $f\in \A^p(\O)$ and by uniqueness of the reproducing function $R_p$, we have $$R_p(z,w)=K(z,w) \text{ and } R_q(z,w)=\overline{R_p(w,z)}=K(z,w),$$ which shows that $K(\cdot,z) \in \A^p(\O)$ for all $z\in\O$.

\erem

\bl \label{Denseness of BK}
 Let $\O$ be a bounded domain and $p,q $ be conjugate exponents of each other. Assume that $\Phi_p, \Phi_q$ are surjective maps. Then $ span\{K(\cdot,z) : z\in \O\}$ is dense in $\A^q(\O)$. 
\el 
\bpf
We will prove the lemma by contradiction. \\
Let $\mathcal{K}=\overline{span\{K(\cdot,z) : z\in \O\}}$, where the closure is with respect to $\A^q$ norm. Suppose that $\mathcal{K}\subsetneq \A^q(\O)$. Then by  Hahn-Banach theorem there exists a linear functional $\varphi \in \A^q(\O)'$ such that $\varphi\neq 0$ and $\varphi_{|\mathcal{K}}\equiv 0$ .\\
By surjectivity of $\Phi_p$, there exists a function $f\in \A^p(\O)$ such that $$\varphi(g)=\<g,f\> for g\in \A^q(\O). $$  By assumption, $\varphi(K(\cdot,z))=0=\<K(\cdot,z),f\>=\overline{f(z)}$. Thus, $f\equiv 0$ contradicting $\phi\neq 0.$
\epf

Next, we consider the following question: If we assume that $\Phi_p$ and $\Phi_q$ are surjective maps, then what can be said about surjectivity of $\Phi_{p'}$ for $q\leq p'\leq p$? An answer to this question lies in the concept of `Complex Interpolation of Banach spaces'.

\section{Complex interpolation of Banach spaces}
We will start this section with a primer on Complex interpolation of Banach Spaces. For more information about Interpolation spaces, cf. 
\cite{Calderón1964}, \cite{cwikel2008lecture}, \cite{bergh2012interpolation}, \cite[\S 10]{zhao2006theory}.

\bedef 
When $X_0,X_1$ are Banach spaces, the ordered pair $(X_0,X_1)$ is called a \textbf{Banach couple} if $X_0$ and $X_1$ are both contained in some Hausdorff topological vector space $\mathscr{X}$ such that the inclusion map $X_i\xhookrightarrow{} \mathscr{X}$ is continuous for $i=0,1$.  \\
We say that a Banach couple $(X_0,X_1)$ is \textbf{regular Banach couple} if $X_0\cap X_1$ is dense in $X_j$ for $j=0,1$. 
\endef 
Let $S$ denote the open strip in the complex plane and $\overline{S}$ denote the closure of $S$: 
$$S=\{x+iy: 0<x<1\} \text{ and } \overline{S}=\{x+iy: 0\leq x\leq1\}.$$
Define the space of continuous functions on $\overline{S}$ to $X_0+X_1$
$$\F(X_0,X_1)=\{f: \overline{S}\rightarrow X_0+X_1: f \text{ is continuous on }\overline{S}, \text{ analytic in }S, f(it)\in X_0, f(1+it)\in X_1\}$$
with the norm $\norm{f}_{\F(X_0,X_1)}=\max\left\{\sup\limits_{t\in \R} \norm{f(it)}_{X_0}, \sup\limits_{t\in \R} \norm{f(1+it)}_{X_1}\right\} <\infty $.

Define $$[X_0,X_1]_\th :=\{f(\th) \in X_0+X_1:  f\in \F(X_0,X_1)\}$$ with the norm $\norm{x}_{[X_0,X_1]_\th}:=\norm{x}_\th= \inf \left\{ \norm{f}_{\F(X_0,X_1)}: f(\th)=x,  f\in \F(X_0,X_1)\right\}$. \\
Then $[X_0, X_1]_\th$ is a Banach space and is called an interpolation space with respect to the Banach couple $(X_0,X_1)$. We note some of the important properties of the interpolation space here:
\begin{itemize}
\item The interpolation space is functorial in the following sense:\\
If $$T: X_0+X_1\rightarrow Y_0+Y_1$$ is a linear operator such that $T$ maps $X_i$ boundedly into $Y_i$ with operator norm less than $C_i$ for $i=0,1$, then $T$ also maps $[X_0,X_1]_\th$ boundedly into $[Y_0,Y_1]_\th$ with norm less than $C_0^{1-\th} C_1^{\th}$.

\item Interpolation spaces also relate to dual spaces in the following way: 
\bl  \cite{Calderón1964}, \cite[Theorem 4.5.1]{bergh2012interpolation}
Let $(X_0,X_1)$ be a regular Banach couple. If $X_0$ and $X_1$ are reflexive spaces, then the dual of the complex interpolation space is obtained by the complex interpolation of the dual spaces, i.e. $$([X_0,X_1]_\th)'= [X_0',X_1']_\th \quad \text{ (with equal norms).}$$
\el 

\item  \cite[Theorem 5.1.1]{bergh2012interpolation} $\L^p$ spaces over a measure space $X$ is a well known example of interpolation spaces. Specifically if  $1\leq p_0< p_1\leq \infty $ and $\frac{1}{p}=\frac{1-\th}{p_0}+\frac{\th}{p_1}$ for some $\th \in(0,1)$, then $$[\L^{p_0}(X), \L^{p_1}(X)]_\th=\L^p(X) \text{ with equal norms.}$$

\end{itemize}

\subsection{Interpolation of Bergman Space}

Throughout this subsection we will assume that $\O$ is  a bounded domain, $p,q$ are conjugate exponents of each other and $1<q<2<p<\infty$. Also assume that $\Phi_p,\Phi_q$ are surjective maps (and hence isomorphisms).

\bt \label{Interpolation}
Let $p_\theta $ be such that $\frac{1}{p_\theta}=\frac{1-\theta}{p}+\frac{\theta}{q}$
for some $\th\in (0,1)$, then 
$$[\A^p(\O), \A^q(\O) ]_{\th}=\A^{p_\th}(\O)$$ with equivalent norms.
\et 
\bpf 
Define the inclusion map $$I : \A^p(\O) \xhookrightarrow{} \L^p(\O )$$
$$I: \A^q(\O) \xhookrightarrow{} \L^q(\O)$$ Then $I$ is bounded linear operator with norm $1$. Further, by the functorial property of complex interpolation, $$I: [\A^p(\O), \A^q(\O)]_\th \rightarrow [\L^p(\O), \L^q(\O)]_\th= \L^{p_\th}(\O).$$\newline 
Thus, $[\A^p(\O), \A^q(\O)]_\th \subset A^{p_\th}(\O)$ and $\norm{g}_{\A^{p_\th}}\leq \norm{g}_{[\th]}$ for all $g\in \A^{p_{\theta}}(\O)$.\\

To prove the other inclusion, define the linear operator $$T: \L^q(\O) \rightarrow \A^q(\O).$$ 
For $g\in \L^q(\O)$, define $Tg  \in \A^q(\O)$ to be the unique holomorphic function such that $$\<f,g\>=\<f, Tg\> \text{ for all } f \in \A^p(\O).$$ Observe that the surjectivity of $\Phi_q$ guarantees the existence of $Tg$ and the injectivity of $\Phi_q$ provides the uniqueness. 
Then $\norm{Tg}_q\lesssim \norm{g}_q$, hence $T$ is a bounded linear operator.\\

Next, we show that $T_{\vert{\L^p(\O)}} \hookrightarrow \A^p(\O) $. Let $g\in \L^p(\O)$, then  $\<\cdot, g\>$ is a linear functional on $\A^q(\O)$ and by  surjectivity of $\Phi_q$, there exists a $g_0 \in \A^p(\O) $ such that $\<h,g_0\>=\<h,g\>$ for all $h\in \A^q(\O)$. Substituting $h=K(\cdot,z)$, we get $Tg=g_0$ i.e. $Tg \in \A^p(\O)$. Therefore,

$$T_{\vert \A^p(\O)}:\L^p(\O)\rightarrow \A^p(\O) \text{ is a bounded operator.}$$
Thus by functorial property of complex interpolation, 
$$ T:[\L^p(\O), \L^q(\O)]_\th= \L^{p_\th} \rightarrow [\A^p(\O), \A^q(\O)]_{\th} $$
which proves $[\A^p(\O), \A^q(\O)]_\th=\A^{p_\th}(\O)$ with $\norm{g}_{p_\theta} \simeq \norm{g}_{[\th]}  $.\epf 

\brem 
The operator T defined above is a projection and thus $\A^p$ and $\A^q$ are complemented subspace of $\L^p$, $\L^q$ respectively. So we can write $\L^p=\A^p\oplus(\A^q)^\perp$ and $\L^q=\A^q \oplus (\A^p)^\perp$, where $(\A^q)^\perp$ and $(\A^p)^\perp$ denotes the annihilators of $\A^q$ and $\A^p$ in $\L^p$ and $\L^q$ respectively.
\erem

\brem \label{duality}
By duality of complex interpolation, $$\A^{q_\th}(\O)'=([\A^q(\O), \A^p(\O)]_\th)'=[\A^p(\O), \A^q(\O)]_\th= \A^{p_\th}(\O)$$
where $\frac{1}{p_\th}=\frac{1-\th}{p}+\frac{\th}{q}$ , $\frac{1}{q_\th}=\frac{1-\th}{q}+\frac{\th}{p}$ and hence $\frac{1}{p_\th}+\frac{1}{q_\th}=1$.
Also, following the proof of duality of complex interpolation shows that $\Phi_{p_\th}$ and $\Phi_{q_\th}$ are also surjective maps i.e.  \textbf{if \bm{$p,q$} are conjugate exponents such that \bm{$\Phi_p,\Phi_q$} are surjective, then for all \bm{$p'$} such that \bm{$q\leq p'\leq p$}, \bm{$\Phi_{p'}$} is a surjective map.}
\erem

\section{Indices of a bounded domain}

For a bounded domain $\O$, B.Y. Chen in \cite{chen2017bergman} defined the \textit{integrabilty index} of the Bergman kernel
$$\b(\O)=\sup \{p\geq 2: K(\cdot,z ) \in \A^p(\O) \text{ for all } z\in \O\}$$
Similarly, Y. Zeytuncu in \cite{zeytuncu2020survey} defined, what we will now call, a \textit{regularity index} of a domain  
$$\mathcal{R}(\O)=\sup\{p\geq 2: \text{ Bergman projection } B: \L^p(\O)\rightarrow\A^p(\O) \text{ is a bounded operator}\}$$
Along similar lines, Remark \ref{duality} now allows us to define the \textit{duality index} of a bounded domain
$$\mathcal{D}(\O)=\sup\{p\geq 2: \Phi_p, \Phi_q  \text{ are surjective maps. }p,q \text{ conjugate of each other}\}$$

\bl
Let $\O$ be a bounded domain then $$\mathcal{D}(\O)\leq \mathcal{R}(\O)\leq \b(\O)$$
\el

\bpf 
Remark \ref{Integrability of BK} says that $\b(\O)\geq \mathcal{D}(\O)$.\\
Also, we know that if $\psi\in C_c^{\infty}(\O)$ is any normalised radially symmetric function with respect to $z\in \O$ then $B(\psi)=K(\cdot,z)$. This shows that $\b(\O)\geq \mathcal{R}(\O).$\\
Now to compare $\mathcal{D}(\O)$ and $\mathcal{R}(\O)$ we will follow the arguments from \cite{Hedenmalm}.\\ Let $p< \mathcal{D}(\O)$ and $q$ be the conjugate exponent of $p$ i.e. $\Phi_p, \Phi_q$ are surjective maps.   Consider the surjective function 
\begin{align*}S : & \L^q(\O) \rightarrow [\A^p(\O)]'  \\
&f \mapsto \<\cdot,f\>
\end{align*}
Then $$ker(S)=\{f\in \L^q(\O) : \<g,f\>=0 \text{ for all } g\in \A^p(\O)\}=[\A^p(\O)]^\perp$$
 where $[\A^p(\O)]^\perp $  represents the annihilators of $\A^p$. By the fundamental theorem of homomorphism, $$[\A^p(\O)]'= A^q(\O) =\L^q(\O)/ [\A^p(\O)]^\perp$$
Therefore $\A^q(\O)$ and $[\A^p(\O)]^\perp$ are complemented subspaces of $\L^q(\O)$
$$\L^q(\O)=\A^q(\O)\oplus [\A^p(\O)]^\perp$$ and there exists a projection $\Pi_q$ ($\Pi_q^2=\Pi_q$) from $\L^q(\O)$ onto $\A^q(\O)$ and  $ker (\Pi_q)=[\A^p(\O)]^\perp$. 
Notice that if $f_1\in \A^2$ and $f_2\in [\A^2]^\perp \subset [\A^p]^\perp$, then for $f=f_1+f_2\in \A^2\oplus [\A^2]^\perp =\L^2(\O) \subset \L^q(\O)$
$$\Pi_q(f_1+f_2)=\Pi_q(f_1)=f_1=B(f_1)=B(f_1+f_2)$$
Thus, $\norm{B(f)}_q\leq \norm{\Pi_q(f)}_q\leq C\norm{f}_q$ on a dense subset of $\L^q(\O)$.
Using the elementary fact that any bounded linear operator defined on a dense subset extends uniquely to a bounded linear operator, we have that $$B:\L^q(\O) \rightarrow\A^q(\O)$$ is bounded linear operator. By duality, it follows that $B:\L^p(\O) \rightarrow\A^p(\O)$ is a bounded linear operator

\epf 

\brem \cite[Theorem 23]{zhao2006theory} states that for a Unit Ball $\B^n$, $\mathcal{D}(\B^n) =\mathcal{R}(\B^n)= \b(\B^n)=\infty$.

\erem 

\subsection{Reinhardt Domains}
  Let $\O\subset \C^n$ be a bounded Reinhardt domain. Then we have that the monomials form a Schauder basis for $\A^p(\O)$. More precisely, denote $S(\O,\A^p)$  to be the set of all $\A^p$ allowable multi-indices, i.e. $$S(\O,\A^p)=\{\a \in \Z^n: e_\a(z):=z_1^{\a_1}\cdots z_n^{\a_n}\in \A^p(\O)\}.$$
Then 

\bt\cite[Theorem 3.11]{chakrabarti2019duality} 
Let $\O$ be a bounded Reinhardt domain, then the set $\{e_\a: \a \in S(\O,\A^p)\}$ forms a Schauder basis for $\A^p(\O)$.
\et 
 Further, as defined in \cite{chakrabarti2019duality}, we call $p\in (1,\infty)$ a \emph{threshold exponent for $\O$} if for for all $\e>0$, $S(\O, \A^p)\neq S(\O, \A^{p+ \e})$ or $S(\O, \A^p)\neq S(\O, \A^{p-\e})$

\bl 
Let $\O\subset \C^n$ be a bounded Reinhardt domain and $2\leq p<\mathcal{D}(\O)$ and q be conjugate to p. Then $S(\O, \A^p)= S(\O, \A^q)=S(\O, \A^2)$.
\el
\bpf
Since $\O$ is a bounded domain, it follows that $S(\O, \A^p)\subset S(\O,\A^2) \subset S(\O,\A^q)$.\\
On contrary, assume $\gamma\in S(\O, \A^q)\setminus S(\O,\A^p)$. Then $$\Phi_p(e_\gamma)(e_\alpha)=\int_\O e_\a \overline{e_\gamma}=0 \text{ for all } \a \in S(\O, \A^p).$$ Hence $e_\gamma\in ker(\Phi_p)$ thereby contradicting  the injectivity of $\Phi_p$ (ref. Lemma \ref{Isomorphic}). 
\epf 
Using the above Lemma,  we observe that there exists an example of a Reinhardt domain $\O$ for which $\mathcal{D}(\O)\neq \mathcal{R}(\O)$. \\\\For $m,n \in \Z^+$, $m,n$ coprime to each other, define the Hartogs triangle 
$$\mathbb{H}_{m/n}=\set{(z_1,z_2) \in \C^2: \abs{z_1}^{m/n}<\abs{z_2}<1 } $$  

Then \cite[Proposition 4.8]{chakrabarti2019duality} states that $2$ is one of the threshold exponents of $\mathbb{H}_{m/n}$. i.e. for all $\e>0$,  $S(\mathbb{H}_{m/n}, \A^2)\neq S(\mathbb{H}_{m/n}, \A^{2-\e})$. Therefore as a consequence of the above Lemma,  $\mathcal{D}(\mathbb{H}_{m/n})=2$.  Further \cite[Theorem 1.1]{edholm2017bergman} states that  $\mathcal{R}(\mathbb{H}_{m/n})=\dfrac{2(m+n)}{m+n-1}$.

\subsection*{Acknowledgements}
The author would like to thank his advisor Prof. Harold Boas for giving valuable advice and feedback during the preparation of this note. He would also like to thank Tanuj Gupta and Siddharth Sabharwal  for useful conversations.

\bibliographystyle{alpha}  
\bibliography{references}

\fontsize{11}{9}\selectfont

\vspace{0.5cm}

\noindent bhatshreedhar33@tamu.edu;

 \vspace{0.2 cm}

\noindent Department of Mathematics, Texas A\&M University, College Station, TX 77843-3368, USA

\end{document}